\documentstyle[amstex,righttag,amscd]{article}
\numberwithin{equation}{section}
\begin{document}
\def\bd{\bold}
\title{On $q$-Clebsch Gordan Rules
and the Spinon Character Formulas
for Affine $C_2^{(1)}$ Algebra}
\author{Yasuhiko Yamada\\
Department of Mathematics\\
Kyushu University}
\maketitle
\begin{abstract}
A $q$-analog of the Clebsch Gordan rules for the tensor products
of the fundamental representations of Yangian is introduced.
Its relation to the crystal base theory and
application to the spinon character formulas are discussed
in case of $C_2^{(1)}$ explicitly. 
\end{abstract}
\section{Introduction}
\par
In 1991, Reshetikhin conjectured the following particle structure
of the one-dimensional $\frak g$-invariant spin chain models
\cite{[R]}.
\begin{itemize}
\item 
There exist $r(={\text{\rm rank}}{\frak g})$ kinds of particles, 
each of them corresponds to a fundamental representation of 
$\frak g$, or more precisely, of Yangian $Y({\frak g})$.
\item
Each particle is a kink. They intertwine two integrable representations
satisfying the admissibility condition, i.e. the restricted fusion rule.
\item
The exchange relations of these particles are described by tensor product
of the Vertex-type of the Yangian and RSOS-type $R$-matrices.
\end{itemize}
The origin of this conjecture is the coincidence of the TBA equations
based on the Bethe equation for the spin chain and on the $S$-matrix
theory.

It is natural to look for a description of the space of states in 
terms of these particles, rather than in terms of symmetry currents.
Such a particle description was first obtained for $sl(2)$ level 1
WZW model, and generalized to higher rank and/or level cases.
Such basis and corresponding character formulas are
called the spinon basis and the spinon character formulas
\cite{[BPS],[BLS],[BS],[ANOT],[NY1],[NY2],[NY3]}.

Generalizing the known cases, we will formulate
spinon character formulas for $C_2^{(1)}$ case. 

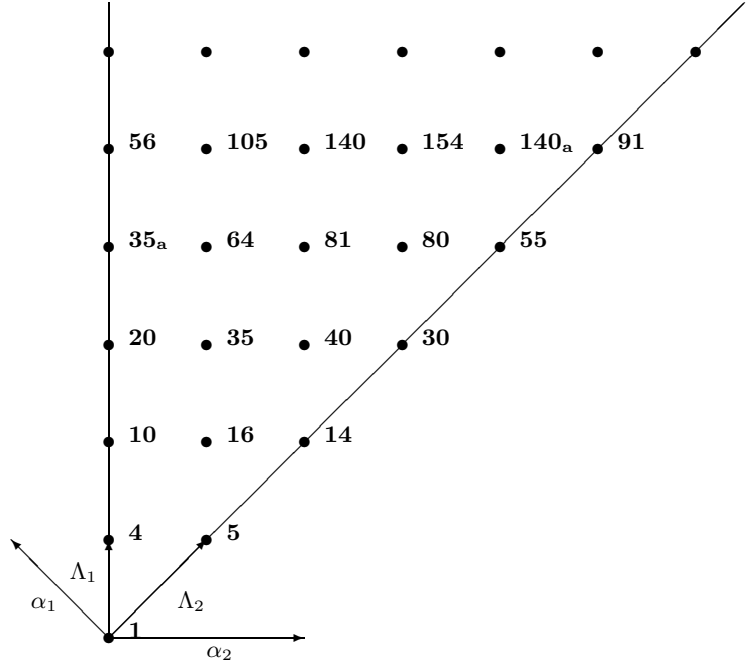
\begin{figure}[h]
\setlength{\unitlength}{1.3mm}
{\small
\begin{picture}(80,90)(-30,-10)
\put(0,0){\line(1,1){65}}
\put(0,0){\line(0,1){65}}
\multiput(0,0)(0,10){7}{\circle*{1}}
\multiput(10,10)(0,10){6}{\circle*{1}}
\multiput(20,20)(0,10){5}{\circle*{1}}
\multiput(30,30)(0,10){4}{\circle*{1}}
\multiput(40,40)(0,10){3}{\circle*{1}}
\multiput(50,50)(0,10){2}{\circle*{1}}
\multiput(60,60)(0,10){1}{\circle*{1}}
\put(0,0){\vector(1,0){20}}
\put(0,0){\vector(1,1){10}}
\put(0,0){\vector(0,1){10}}
\put(0,0){\vector(-1,1){10}}
\put(10,-2){$\alpha_2$}
\put(-8,3){$\alpha_1$}
\put(7,3){$\Lambda_2$}
\put(-4,6){$\Lambda_1$}
\put(2,0){$\bd 1$}
\put(2,10){$\bd 4$}
\put(12,10){$\bd 5$}
\put(2,20){$\bd {10}$}
\put(12,20){$\bd {16}$}
\put(22,20){$\bd {14}$}
\put(2,30){$\bd {20}$}
\put(12,30){$\bd {35}$}
\put(22,30){$\bd {40}$}
\put(32,30){$\bd {30}$}
\put(2,40){$\bd {35}_a$}
\put(12,40){$\bd {64}$}
\put(22,40){$\bd {81}$}
\put(32,40){$\bd {80}$}
\put(42,40){$\bd {55}$}
\put(2,50){$\bd {56}$}
\put(12,50){$\bd {105}$}
\put(22,50){$\bd {140}$}
\put(32,50){$\bd {154}$}
\put(42,50){$\bd {140}_a$}
\put(52,50){$\bd {91}$}
\end{picture}
}
\caption{The representations of $C_2$}
\label{pic:lattice}
\end{figure}

Denote by ${\bd n}$ the irreducible representation of $C_2$ of dimension $n$. 
Some lower dimensional representations are listed in Fig.\ref{pic:lattice}.

Let ${\bd n}^{(k)}$ be the integrable representation of affine
Lie algebra ${C_2^{(1)}}$ of level $k$ such that the $L_0$-ground space
is ${\bd n}$.

For simplicity, let us consider $V={\bd 1}^{(1)}$ case.
The space $V$ has a weight decomposition
$V=\oplus_{n=0}^{\infty} V_n$ w.r.t.  $L_0$-operator. 
And each weight space $V_n$
can be further decomposed into irreducible  ${\frak g}$-modules,
for instance
\begin{eqnarray}
&V_0=&{\bd 1}, \cr
&V_1=&{\bd {10}}, \cr
&V_2=&{\bd 1}+{\bd 5}+{\bd {10}}+{\bd {14}}, \cr
&V_3=&{\bd 1}+{\bd 5}+3 \, {\bd {10}}+{\bd {14}}+{\bd {35}}.
\end{eqnarray}
One can represent this $L_0 \otimes {\frak g}$-module decomposition
formally as 
\begin{eqnarray}
&V=V_0+q V_1+q^2 V_2+\ldots & \cr
&={\bd 1}+q {\bd {10}}+ q^2 ({\bd 1}+{\bd 5}+{\bd {10}}+{\bd {14}})
+\ldots & \cr
&={\bd 1}+{q \over (q)_2} (q {\bd 1}+q {\bd 5}+{\bd {10}})+
{q^2 \over (q)_2} (q^2 {\bd 1}+q {\bd {10}}+{\bd {14}})+\ldots,&
\end{eqnarray}
where $(q)_n=(1-q)(1-q^2)\cdots (1-q^n)$, and the expression
$1/(q)_n$ should be considered as a formal power series in $q$.

The last expression can be interpreted in terms of the particle
picture.
The first term ${\bd 1}$ is the contribution of the unique ground state.
The second term is the contribution of the two particles each of them
in ${\bd 4}$ representation, and the third term comes from ${\bd 5}^2$.
One can recognize the expression
\begin{eqnarray}
&[{\bd 4}^2]=q {\bd 1}+q {\bd 5}+{\bd {10}},& \cr
&[{\bd 5}^2]=q^2 {\bd 1}+q {\bd {10}}+{\bd {14}},&
\end{eqnarray}
are $q$-analog of the Clebsch Gordan rules.

For any integrable representation $V$ of $C_2^{(1)}$,
it is natural to expect the following formula
\begin{equation}
V=\sum_{n,m=0}^{\infty} {P(q) \over (q)_n (q)_m} 
[{\bd 4}^n \otimes {\bd 5}^m],
\end{equation}
where $P(q)$ is some polynomial.

The aim of this note is to formulate this expectation
as explicit as possible, and check them.
The paper is organized as follows.
In section 2, the first definition of the $q$-Clebsch Gordan rules
are formulated using a result from the Bethe ansatz
\cite{[KR],[KR1],[KNS],[Kir]}.
In section 3, we prepare some basic data from the crystal theory
in $C_2^{(1)}$ case. In section 4, we propose the second definition
of the $q$-Clebsch Gordan rules including the restricted case,
and apply them to the character formulas. Section 5 is for
summary and discussions.

\section{$q$-Clebsch Gordan rules via Bethe ansatz}
\par
Let $\frak g$ be a finite dimensional simple Lie algebra of rank $r$.
Let $\Lambda_{a}$, $\alpha_{a}$ and $\alpha_a^{\vee}=2 \alpha_a/\alpha_a^2$, 
($a=1,\cdots,r$) be the fundamental weighs, simple roots and co-roots.
We normalize them so that $\theta^2=2$ for the long root $\theta$.

For a pair of dominant integral weights $\lambda, \mu$,
we define a polynomial $M_{\lambda, \mu}(q)$ 
by the following formula,
\begin{equation}
M_{\lambda,\mu}(q)=\sum_{m} q^{c(m)}
\prod_{a=1}^{r} \prod_{i=1}^{\infty}
\left[ \begin{array}{c}
P^{a}_i (m)+m^{a}_i \\
m^{a}_i
\end{array} \right]_q .
\end{equation}  
The sum is taken over the nonnegative integers 
$m^{a}_i$ ($ a=1, \cdots, r$ and $i=1,2,\cdots$) such that
\begin{equation}
\lambda=\mu-\sum_{a=1}^{r}
(\sum_{i=1}^{\infty} i m^{a}_i) \alpha_a.
\end{equation}
The number $P^{a}_i(m)$ is defined by
\begin{equation}
P^{a}_i (m)=(\alpha_a^{\vee},\mu)-\sum_{b=1}^{r} \sum_{j=1}^{\infty}
\Phi^{a,b}_{i,j} m^{b}_j,
\end{equation}
\begin{equation}
\Phi^{a,b}_{i,j}=2 {(\alpha_a, \alpha_b) \over \alpha_a^2 \alpha_b^2}
\min\{i \alpha_a^2, j \alpha_b^2\},
\end{equation}
The $q$-binomial coefficients are defined as
\begin{equation}
\left[ \begin{array}{c} a\\ b \end{array} \right]_q=
{(q)_{a+b} \over (q)_a (q)_b},
\end{equation}
with $(q)_a=(1-q)(1-q^2)\cdots(1-q^a)$, 
and finally, the charge $c$ is given by
\begin{equation}
c(m)={1 \over 2} \sum_{a,b=1}^{r} \sum_{i,j=1}^{\infty} 
\Phi^{a,b}_{i,j} m^{a}_i m^{b}_j.
\end{equation}

The definition of $M_{\lambda,\mu}$ comes from the analysis 
of the Bethe ansatz equation \cite{[KR],[KR1],[KNS]}. 
We denote by $V(\lambda)$ the irreducible highest weight
representation of $\frak g$ with highest weight $\lambda$.
Denote by $W_a$ the minimal irreducible representation of 
$U_q({\frak g})$ or Yangian $Y({\frak g})$ such that 
$W_a \supseteq V(\Lambda_a)$.
Consider irreducible decomposition of tensor products
of these representations $W_a$ as ${\frak g}$-modules
\begin{equation}
W_1^{n_1} \otimes W_2^{n_2} \otimes \cdots \otimes
W_r^{n_r}=\sum_{\lambda} M_{\lambda, \mu} V(\lambda),
\end{equation}
where $n_a=(\alpha_a^{\vee},\mu)$.

By the analysis of the Bethe ansatz equation,
the following formula for the multiplicity $M_{\lambda,\mu}$
is obtained \cite{[KR1]}(see also \cite{[KNS]})
\begin{equation}
M_{\lambda, \mu}=\sum_{m} \prod_{a=1}^{r} \prod_{i=1}^{\infty}
\left[ \begin{array}{c}
P^{a}_i (m)+m^{a}_i \\
m^{a}_i
\end{array} \right]. 
\end{equation}
The above definition of $M_{\lambda,\mu}(q)$ 
is a $q$-analogue of the formula $M_{\lambda,\mu}$.

In this sense, we define the $q$-Clebsch Gordan rules as follows
\begin{equation}
[W_1^{n_1} \otimes W_2^{n_2} \otimes \cdots \otimes
W_r^{n_r}]=\sum_{\lambda} M_{\lambda, \mu}(q) V(\lambda),
\end{equation}

In case of ${\frak g}={\frak sl}_{r+1}$,  
$M_{\lambda,\mu}(q)$ is nothing but the Kostka polynomial
$K_{\lambda, \mu'}(q)$, where $\mu'$ is the transposition of $\mu$ 
as a Young tableaux \cite{[M],[Kir]}.
And in this case, a crystal theoretic formula of the Kostka polynomial 
is known \cite{[NY2],[D],[DF],[KMOTU]}.
We will generalize it for ${\frak g}=C_2$ case.
\section{The Isomorphism and the Energy function}

In what follows, we concentrate on the $C_2^{(1)}$ case.
We call a sequence '$p_1 p_2 \ldots p_m$'
of letters $p_i \in \{a,b,c,d,e,1,2,3,4\}$ as a word.
In the crystal theory, ${\bd 4}=\{1,2,3,4\}$ and
${\bd 5}=\{a,b,c,d,e\}$ are the crystal base for the two 
fundamental representations of $C_2$ (and also for $C_2^{(1)}$) 
\cite{[K1],[KMN1],[KMN2],[KN]}.
Each word labels the crystal base of corresponding tensor product.
For instance, '$abc12$' is a base in ${\bd 5}^3 \otimes {\bd 4}^2$.

By considering a word as a label of the crystal base,
one can define an isomorphism and a weight function on them \cite{[NY2]}.
Here we give their explicit description for the $C_2$ case.
The relevant crystal graphs are given in Appendix 1.

The isomorphism is a rule to identify words corresponding to
the equivalent vertices on the crystal graphs.
The rule consists of local identification rules concerning
two adjacent letters.

For ${\bd 4}^2$ and ${\bd 5}^2$ cases, the isomorphism is trivial
(identity) map.

The isomorphism ${\bd 5}\otimes{\bd 4}={\bd 4}\otimes{\bd 5}$
is given by
\begin{table}[h]
\begin{center}
\begin{tabular}{ccccc}
$a1=1a$ & $b1=1b$ & $c1=2b$ & $d1=2c$ & $e1=3c$ \\
$a2=2a$ & $b2=3a$ & $c2=4a$ & $d2=2d$ & $e2=3d$ \\
$a3=1c$ & $b3=3b$ & $c3=4b$ & $d3=4c$ & $e3=3e$ \\
$a4=1d$ & $b4=1e$ & $c4=2e$ & $d4=4d$ & $e4=4e$
\end{tabular}
\end{center}
\end{table}

The Energy function $H$ is a ${\bd Z}$-valued
function defined for two letters.

Explicitly, they are given as follows.
\begin{equation}
H(i,j)=
\begin{cases}
0, & \text{if $ij=11,21,22,31,32,33,41,42,43,44$} \\
1, & \text{if $ij=12,13,23,24,34,14$.}
\end{cases}
\end{equation}
\begin{equation}
H(u ,v)= 
\begin{cases}
0, & \text{if $uv=aa,ba,bb,ca,cb,da,db,dc,dd,ea,eb,ec,ed,ee$} \\
1, & \text{if $uv=ab,ac,ad,bc,bd,be,cc,cd,ce,de$} \\
2, & \text{if $uv=ae$.} 
\end{cases}
\end{equation}
\begin{equation}
H(u,j)=
\begin{cases}
0, & \text{if $uj=a1,a2,b1,b2,b3,c1,c2,c3,d1,d2,d3,d4,e1,e2,e3,e4$} \\
1, & \text{if $uj=a3,a4,b4,c4$.}
\end{cases}  
\end{equation}
\begin{equation}
H(j,u)=
\begin{cases}
0, & \text{if $ju=1a,1b,2a,2b,2c,2d,3a,3b,3c,3d,3e,4a,4b,4c,4d,4e$} \\
1, & \text{if $ju=1c,1d,1e,2e$.}
\end{cases}  
\end{equation}
Here, the sign of $H$ is different from the conventional one.

The weight of a word '$p_1 p_2 \ldots p_m$' is defined as the sum of
the index, ${\text{\rm ind}}(p_i)$, of each letter $p_i$.
\begin{equation}
{\text{\rm weight}}(p_1 p_2 \ldots p_m)=
\sum_{i=1}^{m} {\text{\rm ind}}(p_i).
\end{equation}
${\text{\rm ind}}(p_i)$ is defined as follows.
It is easy to describe it by an example.

Let us consider the index of the 4-th letter '2' in the word '$ae123$'.
By using the isomorphism, one can move the letter
'2' to the top of the sequence as follows
\begin{equation}
ae(1{\underline 2})3 \rightarrow 
a(e{\underline 1})23 \rightarrow 
(a{\underline 3})c23 \rightarrow 
{\underline 1}cc23.
\end{equation}
The index is defined as the sum of the Energy function $H$
corresponding to the local exchange at each step
\begin{equation}
H(1,2)+H(e,1)+H(a,3)=1+0+1=2.
\end{equation}
By similar counting, the index of each letter in the word
'$ae123$' is determined as
\begin{equation}
{\text{\rm ind}}(a)=0, \quad
{\text{\rm ind}}(e)=2, \quad
{\text{\rm ind}}(1)=1, \quad
{\text{\rm ind}}(2)=2, \quad
{\text{\rm ind}}(3)=3.
\end{equation}
Hence, the weight of the word '$ae123$' is $0+2+1+2+3=8$.

\section{$q$-Clebsch Gordan rules via Crystal theory}
\par

Here we propose another definition of the $q$-Clebsch Gordan rules.
A path is a sequence of arrows on the weight diagram.
Each arrow corresponds to one word '$1,2,3,4,a,b,c,d$' or '$e$' 
as depicted in Fig.2.

Let us recall the (level $k$-restricted) fusion rules.
For $k=1$ case, there are three integrable representation of $C_2^{(1)}$, 
denoted by ${\bd 1}^{(1)}$, ${\bd 4}^{(1)}$ and ${\bd 5}^{(1)}$ 
respectively.
Their fusion rules are,
\begin{equation}
{\bd 1}^{(1)} \overset{1}{\rightarrow} {\bd 4}^{(1)},\quad
{\bd 4}^{(1)} \overset{2}{\rightarrow} {\bd 5}^{(1)},\quad
{\bd 4}^{(1)} \overset{4}{\rightarrow} {\bd 1}^{(1)},\quad
{\bd 5}^{(1)} \overset{3}{\rightarrow} {\bd 4}^{(1)},
\end{equation}
for ${\bd 4}$ and
\begin{equation}
{\bd 1}^{(1)} \overset{a}{\rightarrow} {\bd 5}^{(1)},\quad
{\bd 4}^{(1)} \overset{c}{\rightarrow} {\bd 4}^{(1)},\quad
{\bd 5}^{(1)} \overset{e}{\rightarrow} {\bd 1}^{(1)},
\end{equation}
for ${\bd 5}$.

One can express the fusion rule by the diagram Fig.\ref{fig:fusion1}.
For the higher level $k$ cases, one can represent the 
fusion rule by a similar diagrams.
As an example, we present the diagram for $k=3$ (Fig.\ref{fig:fusion3}).
If the level $k$ is large enough, the fusion rule is noting but the
ordinary Clebsch Gordan rule (i.e. the unrestricted fusion rule).

A path that satisfy the restricted [unrestricted] fusion rule
is called the restricted [unrestricted] fusion path.
A weight of a path is defined as the weight of corresponding word. 

We usually put the starting point of the path on ${\bd 1}^{(k)}$. 
If the end point is ${\bd n}^{(k)}$, we call the path a ${\bd n}^{(k)}$-path.

For instance, there exist two level $1$-restricted ${\bd 4}^{(1)}$-paths
in ${\bd 5}^2 \otimes {\bd 4}^3$, such as
\begin{eqnarray}
&ae141=&[{\bd 1}^{(1)} 
\overset{a}{\rightarrow} {\bd 5}^{(1)}
\overset{e}{\rightarrow} {\bd 1}^{(1)}
\overset{1}{\rightarrow} {\bd 4}^{(1)}
\overset{4}{\rightarrow} {\bd 1}^{(1)}
\overset{1}{\rightarrow} {\bd 4}^{(1)}], \cr
&ae123=&[{\bd 1}^{(1)} 
\overset{a}{\rightarrow} {\bd 5}^{(1)}
\overset{e}{\rightarrow} {\bd 1}^{(1)}
\overset{1}{\rightarrow} {\bd 4}^{(1)}
\overset{2}{\rightarrow} {\bd 5}^{(1)}
\overset{3}{\rightarrow} {\bd 4}^{(1)}].
\end{eqnarray}

Now we define the level $k$-restricted $q$-Clebsch Gordan rules as,
\begin{equation}
[{\bd 4}^n \otimes {\bd 5}^m]^{(k)}=
\sum_{\bd n} \sum_{p} q^{{\text {\rm weight}}(p)} {\bd n}^{(k)},
\end{equation}
where the sum is taken over all the restricted ${\bd n}^{(k)}$-path
in the tensor product ${\bd 4}^n \otimes {\bd 5}^m$.

By the isomorphism described in previous section, there are
different (but equivalent) way to represent tensor products.
For instance, one has
$
{\bd 5} \otimes {\bd 5} \otimes {\bd 4} \otimes {\bd 4} \otimes {\bd 4}
\rightarrow
{\bd 5} \otimes {\bd 4} \otimes {\bd 4} \otimes {\bd 5} \otimes {\bd 4}
$
which maps the two paths above as
\begin{equation}
ae141 \rightarrow a32e1, \quad
ae123 \rightarrow a34a3.
\end{equation}
Since the $H$-function is invariant w.r.t. the isomorphism,
the $q$-Clebsch Gordan rules are independent of this choice, i.e.
$
[{\bd 5} \otimes {\bd 5} \otimes {\bd 4} \otimes {\bd 4} \otimes {\bd 4}]^{(k)}
=
[{\bd 5} \otimes {\bd 4} \otimes {\bd 4} \otimes {\bd 5} \otimes {\bd 4}]^{(k)}
$.

We list some examples of level $1$-restricted $q$-Clebsch Gordan rules.

2-particles
\begin{eqnarray}
[{\bd 4}^2]^{(1)}&=&q\,{\bd 1}^{(1)} + q\,{\bd 5}^{(1)}, \cr
[{\bd 4} \otimes {\bd 5}]^{(1)}&=&q\,{\bd 4}^{(1)}, \cr
[{\bd 5}^2]^{(1)}&=&{q^2}\,{\bd 1}^{(1)},
\end{eqnarray}

3-particles
\begin{eqnarray}
[{\bd 4}^3]^{(1)}&=& 
\left({q^2} + {q^3} \right) \,{\bd 4}^{(1)}, \cr
[{\bd 4}^2 \otimes {\bd 5}]^{(1)}&=&
{q^3}\,{\bd 1}^{(1)} + {q^2} \,{\bd 5}^{(1)}, \cr
[{\bd 4} \otimes {\bd 5}^2]^{(1)}&=&{q^3} \,{\bd 4}^{(1)}, \cr
[{\bd 5}^3]^{(1)}&=&{q^4} \,{\bd 5}^{(1)},
\end{eqnarray}

4-particles
\begin{eqnarray}
[{\bd 4}^4]^{(1)}&=&
\left( {q^4} + {q^6} \right) \,{\bd 1}^{(1)} + 
\left( {q^4} + {q^5} \right) \,{\bd 5}^{(1)}, \cr
[{\bd 4}^3 \otimes {\bd 5}]^{(1)}&=&
\left( {q^4} + {q^5} \right) \,{\bd 4}^{(1)}, \cr
[{\bd 4}^2 \otimes {\bd 5}^2]^{(1)}&=&
{q^5} \,{\bd 1}^{(1)} + {q^5} \,{\bd 5}^{(1)}, \cr
[{\bd 4} \otimes {\bd 5}^3]^{(1)}&=&{q^6} \,{\bd 4}^{(1)}, \cr
[{\bd 5}^4]^{(1)}&=&{q^8} \,{\bd 1}^{(1)}.
\end{eqnarray}

For the unrestricted case, the $q$-Clebsch Gordan rules are listed in the
appendix 3. From these explicit examples, we observe the following
\vskip2mm
\noindent
{\bf Conjecture 1.}

The unrestricted $q$-Clebsch Gordan rules are coincide with the 
$q$-Clebsch Gordan rules defined in section 2.
\vskip2mm

\begin{figure}[h]
\setlength{\unitlength}{1.3mm}
{\small
\begin{picture}(120,30)(10,-10)
\put(10,10){\circle*{1}}
\put(30,10){\circle*{1}}
\put(20,0){\circle*{1}}
\put(20,20){\circle*{1}}
\put(21,10){\vector(1,0){8}}
\put(20,11){\vector(0,1){8}}
\put(19,10){\vector(-1,0){8}}
\put(20,9){\vector(0,-1){8}}
\put(21,15){1}
\put(21,5){4}
\put(15,11){3}
\put(25,11){2}
\put(20,-5){${\bd 4}$}
\put(70,10){\circle*{1}}
\put(60,0){\circle*{1}}
\put(60,20){\circle*{1}}
\put(80,0){\circle*{1}}
\put(80,20){\circle*{1}}
\put(71,11){\vector(1,1){8}}
\put(71,9){\vector(1,-1){8}}
\put(69,11){\vector(-1,1){8}}
\put(69,9){\vector(-1,-1){8}}
\put(66,10){\oval(8,4)}
\put(65,8){\vector(1,0){2}}
\put(76,14){a}
\put(76,5){d}
\put(66,15){b}
\put(66,4){e}
\put(60,10){c}
\put(70,-5){${\bd 5}$}
\end{picture}
\caption{The weights and corresponding arrows}
\label{fig:arrow}
}
%
\setlength{\unitlength}{1.3mm}
{\small
\begin{picture}(120,30)(0,-10)
\put(0,0){\circle*{1}}
\put(0,10){\circle*{1}}
\put(10,10){\circle*{1}}
\put(0,0){\line(0,1){10}}
\put(0,10){\line(1,0){10}}
\put(2,3){$\bd 1$}
\put(2,13){$\bd 4$}
\put(12,13){$\bd 5$}
\put(5,-5){$\bd 4$}
\put(50,0){\circle*{1}}
\put(50,10){\circle*{1}}
\put(60,10){\circle*{1}}
\put(50,0){\line(1,1){10}}
\put(48,10){\oval(4,2)}
\put(50,3){$\bd 1$}
\put(50,13){$\bd 4$}
\put(60,13){$\bd 5$}
\put(55,-5){$\bd 5$}
\end{picture}
}
\caption{The fusion rule for level $k=1$}
\label{fig:fusion1}
%
\setlength{\unitlength}{1.3mm}
{\small
\begin{picture}(120,50)(0,-10)
\put(0,0){\line(0,1){30}}
\put(10,10){\line(0,1){20}}
\put(20,20){\line(0,1){10}}
\put(0,10){\line(1,0){10}}
\put(0,20){\line(1,0){20}}
\put(0,30){\line(1,0){30}}
\multiput(0,0)(0,10){4}{\circle*{1}}
\multiput(10,10)(0,10){3}{\circle*{1}}
\multiput(20,20)(0,10){2}{\circle*{1}}
\multiput(30,30)(0,10){1}{\circle*{1}}
\put(2,3){$\bd 1$}
\put(2,13){$\bd 4$}
\put(12,13){$\bd 5$}
\put(2,23){$\bd {10}$}
\put(12,23){$\bd {16}$}
\put(22,23){$\bd {14}$}
\put(2,33){$\bd {20}$}
\put(12,33){$\bd {35}$}
\put(22,33){$\bd {40}$}
\put(32,33){$\bd {30}$}
\put(15,-5){$\bd 4$}
\put(50,0){\line(1,1){30}}
\put(50,10){\line(1,1){20}}
\put(50,20){\line(1,1){10}}
\put(60,10){\line(-1,1){10}}
\put(60,20){\line(-1,1){10}}
\put(70,20){\line(-1,1){10}}
\put(48,10){\oval(4,2)}
\put(48,20){\oval(4,2)}
\put(48,30){\oval(4,2)}
\put(58,20){\oval(4,2)}
\put(58,30){\oval(4,2)}
\put(68,30){\oval(4,2)}
\multiput(50,0)(0,10){4}{\circle*{1}}
\multiput(60,10)(0,10){3}{\circle*{1}}
\multiput(70,20)(0,10){2}{\circle*{1}}
\multiput(80,30)(0,10){1}{\circle*{1}}
\put(50,3){$\bd 1$}
\put(50,13){$\bd 4$}
\put(60,13){$\bd 5$}
\put(50,23){$\bd {10}$}
\put(60,23){$\bd {16}$}
\put(70,23){$\bd {14}$}
\put(50,33){$\bd {20}$}
\put(60,33){$\bd {35}$}
\put(70,33){$\bd {40}$}
\put(80,33){$\bd {30}$}
\put(65,-5){$\bd 5$}
\end{picture}
}
\caption{The fusion rule for level $k=3$}
\label{fig:fusion3}
\end{figure}

Finally, we will formulate the spinon character formulas
as an application of the $q$-Clebsch Gordan rules.
As in the section 1, let ${\bd n}^{(k)}$ be the integrable
representation of affine $C_2^{(1)}$ algebra with level $k$
and with $L_0$-ground space ${\bd n}$.
Using the (restricted) $q$-Clebsch Gordan rules, one can formulate
the particle structure of the space ${\bd n}^{(k)}$ as follows,
\vskip2mm
\noindent
{\bf Conjecture 2.}

The representation ${\bd n}^{(k)}$ has the following particle
decomposition 
\begin{equation}
{\bd n}^{(k)}=\sum_{n,m=0}^{\infty}
{P(q) \over (q)_n (q)_m} [{\bd 4}^n \otimes {\bd 5}^m],
\end{equation}
where $P(q)$ is the coefficient of the ${\bd n}^{(k)}$
in the restricted $q$-Clebsch Gordan rule 
$[{\bd 4}^n \otimes {\bd 5}^m]^{(k)}$.
\vskip2mm
This formula is consistent with the known 
results for ${\frak g}={\frak {sl}}_{r+1}$ \cite{[BS],[NY3]}.

\section{Summary and Discussions}

We proposed two definitions of the $q$-Clebsch Gordan rules based
on the Bethe ansatz (section 2) and crystal theory (section 4), 
and apply them to the character formulas for affine Lie algebra $C_2^{(1)}$.
These two definitions are expected to be equivalent.
 
Similar results has been proved in case of ${\frak g}={\frak sl_{r+1}}$
\cite{[NY2]} by using known properties on the Kostka polynomials.
We remark a crucial difference between ${\frak sl_{r+1}}$ and
other cases.
There exist a generalization of the Kostka polynomials
for all Lie algebras as the $q$-analog of the 'weight multiplicity'.
However, as we clarified in $C_2^{(1)}$ case explicitly,
what we need for the character formula is the 
$q$-analog of the 'tensor product multiplicity'.
The latter one is different from the former except for 
${\frak g}={\frak sl_{r+1}}$ cases.

Another interesting problem is to find the Bethe ansatz formula
for the restricted fusion path counting.

\vskip 1cm

{\it Acknowledgment.}
I would like to thank A. Nakayashiki and A. N. Kirillov
for valuable discussion.

\vfill \break

\centerline{Appendix 1 :  Crystal graphs}
\vskip1cm
The crystal graphs for the fundamental representations 
${\bd 4}$ and ${\bd 5}$ are
\begin{equation*}
\begin{CD}
1 @>1>> 2 @>2>> 3 @>1>> 4.
\end{CD}
\end{equation*}
\begin{equation*}
\begin{CD}
a @>2>> b @>1>> c @>1>> d @>2>> e.
\end{CD}
\end{equation*}
For the tensor product 
${\bd 4}\otimes{\bd 4}={\bd 1}\oplus{\bd 5}\oplus{\bd {10}}$, one has
\begin{equation*}
\begin{CD}
11 @>1>> 21 @>2>> 31 @>1>> 41   \\
@.     @V1VV      @.      @V1VV \\
12 @.    22 @>2>> 32  @.   42   \\
@V2VV    @.      @V2VV    @V2VV \\
13 @>1>> 23 @.    33 @>1>> 43   \\
@.     @V1VV      @.      @V1VV \\
14 @.    24 @>2>> 34 @.    44.
\end{CD}
\end{equation*}
Similarly, for 
${\bd 5}\otimes{\bd 5}={\bd 1}\oplus{\bd {10}}\oplus{\bd {14}}$,
\begin{equation*}
\begin{CD}
aa @>2>> ba @>1>> ca @>1>> da @>2>> ea   \\
@.      @V2VV     @.       @.      @V2VV \\
ab @.    bb @>1>> cb @>1>> db @.    eb   \\
@V1VV    @.       @.      @V1VV    @V1VV \\
ac @>2>> bc @>1>> cc @.    dc @.    ec   \\
@V1VV    @.      @V1VV    @V1VV    @V1VV \\
ad @>2>> bd @.    cd @.    dd @>2>> ed   \\
@.      @V2VV    @V2VV     @.      @V2VV \\
ae @.    be @>1>> ce @>1>> de @.    ee.
\end{CD}
\end{equation*}
For the tensor products
${\bd 4}\otimes{\bd 5}={\bd 5}\otimes{\bd 4}={\bd {4}}\oplus{\bd {16}}$,
one has
\begin{equation*}
\begin{CD}
1a @>1>> 2a @>2>> 3a @>1>> 4a   \\
@V2VV    @.      @V2VV    @V2VV \\
1b @>1>> 2b @.    3b @>1>> 4b   \\
@.      @V1VV     @.      @V1VV \\
1c @.    2c @>2>> 3c @.    4c   \\
@V1VV   @V1VV    @V1VV    @V1VV \\
1d @.    2d @>2>> 3d @.    4d   \\
@V2VV    @.      @V2VV    @V2VV \\
1e @>1>> 2e @.    3e @>1>> 4e.
\end{CD}
\end{equation*}
and
\begin{equation*}
\begin{CD}
a1 @>2>> b1 @>1>> c1 @>1>> d1 @>2>> e1   \\
@V1VV    @.       @.      @V1VV    @V1VV \\
a2 @>2>> b2 @>1>> c2 @.    d2 @>2>> e2   \\
@.      @V2VV    @V2VV     @.      @V2VV \\
a3 @.    b3 @>1>> c3 @>1>> d3 @.    e3   \\
@V1VV    @.       @.      @V1VV    @V1VV \\
a4 @>2>> b4 @>1>> c4 @.    d4 @>2>> e4.
\end{CD}
\end{equation*}

\vfill \break
\centerline{Appendix 2 :  Restricted fusion path}
\vskip1cm
\begin{table}[h]
\caption{level $k=1$ ${\bd 1}^{(1)}$-path}
\begin{center}
\begin{tabular}{|c|cccc|}
\noalign{\hrule height1.2pt}
particle content && path(weight) && \\
\noalign{\hrule height1.2pt}
${\bd 1}$ & $\phi_{0}$ &&& \\
\hline
${\bd 4}^2$ & $14_1$ &&& \\
${\bd 5}^2$ & $ae_2$ &&& \\
\hline
${\bd 4}^2\otimes{\bd 5}$ & $a34_3$ &&& \\
\hline
${\bd 4}^4$ & $1414_4$ & $1234_6$ && \\
${\bd 4}^2\otimes{\bd 5}^2$ & $ae14_5$ &&& \\
${\bd 5}^4$ & $aeae_8$ &&& \\
\hline
${\bd 4}^4\otimes{\bd 5}$ & $a3234_7$ & $a3414_8$ && \\
${\bd 4}^2\otimes{\bd 5}^3$ & $aea34_9$ &&& \\
\hline
${\bd 4}^6 $ & $141414_9$
             & $141234_{11}$
             & $123234_{12}$
             & $123414_{13}$ \\
${\bd 4}^4\otimes{\bd 5}^2$ & $ae1414_{10}$
                            & $ae1234_{12}$ && \\
${\bd 4}^2\otimes{\bd 5}^4$ & $aeae14_{13}$ &&& \\
${\bd 5}^6$ & $aeaeae_{18}$ &&& \\
\noalign{\hrule height1.2pt}
\end{tabular}
\end{center}
\end{table}

\begin{table}[h]
\caption{level $k=1$ ${\bd 4}^{(1)}$-path}
\begin{center}
\begin{tabular}{|c|cccc|}
\noalign{\hrule height1.2pt}
particle content && path(weight) && \\
\noalign{\hrule height1.2pt}
${\bd 4}$ & $1_{0}$ &&& \\
\hline
${\bd 4}\otimes{\bd 5}$ & $a3_1$ &&& \\
\hline
${\bd 4}^3$ & $141_2$ & $123_3$ && \\
${\bd 4}\otimes{\bd 5}^2$ & $ae1_3$ &&& \\
\hline
${\bd 4}^3\otimes{\bd 5}$ & $a323_4$ & $a341_5$ && \\
${\bd 4}\otimes{\bd 5}^3$ & $aea3_6$ &&& \\
\hline
${\bd 4}^5$ & $14141_6$ & $14123_7$ & $12323_8$ & $12341_9$ \\
${\bd 4}^3\otimes{\bd 5}^2$ & $ae141_7$ & $ae123_8$ && \\
${\bd 4}\otimes{\bd 5}^4$ & $aeae1_{10}$ &&& \\
\hline   
${\bd 4}^5\otimes{\bd 5}$ & $a32323_9$
             & $a32341_{10}$
             & $a34141_{11}$
             & $a34123_{12}$ \\
${\bd 4}^3\otimes{\bd 5}^3$ & $aea323_{11}$
                            & $aea341_{12}$ && \\
${\bd 4}\otimes{\bd 5}^5$ & $aeaea3_{15}$ &&& \\
\noalign{\hrule height1.2pt} 
\end{tabular}
\end{center}
\end{table}

\vfill \break

\begin{table}[h]
\caption{level $k=1$ ${\bd 5}^{(1)}$-path}
\begin{center}
\begin{tabular}{|c|cccc|}
\noalign{\hrule height1.2pt}
particle content && path(weight) && \\
\noalign{\hrule height1.2pt}
${\bd 5}$ & $a_{0}$ &&& \\
\hline
${\bd 4}^2$ & $12_1$ &&& \\
\hline
${\bd 4}^2\otimes{\bd 5}$ & $a32_2$ &&& \\
${\bd 5}^3$ & $aea_4$ &&& \\
\hline
${\bd 4}^4$ & $1412_4$ & $1232_5$ && \\
${\bd 4}^2\otimes{\bd 5}^2$ & $ae12_5$ &&& \\
\hline
${\bd 4}^4\otimes{\bd 5}$ & $a3232_6$ & $a3412_8$ && \\
${\bd 4}^2\otimes{\bd 5}^3$ & $aea32_8$ &&& \\
${\bd 5}^5$ & $aeaea_{12}$ &&& \\
\hline   
${\bd 4}^6 $ & $141412_9$
             & $141232_{10}$
             & $123232_{11}$
             & $123412_{13}$ \\
${\bd 4}^4\otimes{\bd 5}^2$ & $ae1412_{10}$
                            & $ae1232_{11}$ && \\
${\bd 4}^2\otimes{\bd 5}^4$ & $aeae12_{13}$ &&& \\
\noalign{\hrule height1.2pt} 
\end{tabular}
\end{center}
\end{table}

\vfill \break

\centerline{Appendix 3 :  $q$-Clebsch Gordan rules}
\vskip 1cm

2-particles
\begin{eqnarray*}
[{\bd 4}^2]&=&q\,{\bd 1} + q\,{\bd 5} + {\bd {10}}, \cr
[{\bd 4} \otimes {\bd 5}]&=&q\,{\bd 4} + {\bd {16}}, \cr
[{\bd 5}^2]&=&{q^2}\,{\bd 1} + q\,{\bd {10}} + {\bd {14}},
\end{eqnarray*}

3-particles
\begin{eqnarray*}
[{\bd 4}^3]&=& 
\left( q + {q^2} + {q^3} \right) \,{\bd 4} + 
\left( q + {q^2} \right) \,{\bd {16}} + {\bd {20}}, \cr
[{\bd 4}^2 \otimes {\bd 5}]&=&
{q^3}\,{\bd 1} + \left( q + {q^2} \right) \,{\bd 5} + 
\left( q + {q^2} \right) \,{\bd {10}} + q\,{\bd {14}} + {\bd {35}}, \cr
[{\bd 4} \otimes {\bd 5}^2]&=&
\left( {q^2} + {q^3} \right) \,{\bd 4} + 
\left( q + {q^2} \right) \,{\bd {16}} + 
q\,{\bd {20}} + {\bd {40}}, \cr
[{\bd 5}^3]&=&
\left( {q^2} + {q^3} + {q^4} \right) \,{\bd 5} + {q^3}\,{\bd {10}} + 
\left( q + {q^2} \right) \,{\bd {35}} + {\bd {30}}
\end{eqnarray*}

4-particles
\begin{eqnarray*}
[{\bd 4}^4]&=&
\left( {q^2} + {q^4} + {q^6} \right) \,{\bd 1} + 
\left( {q^2} + {q^3} + 2\,{q^4} + {q^5} \right) \,{\bd 5} + \cr &&
\left( q + {q^2} + 2\,{q^3} + {q^4} + {q^5} \right) \,{\bd {10}} +  
\left( {q^2} + {q^4} \right) \,{\bd {14}} + \cr &&
\left( q + {q^2} + {q^3} \right) \,{\bd {35}} + {\bd {35}_{a}}, \cr
[{\bd 4}^3 \otimes {\bd 5}]&=&
\left( {q^2} + {q^3} + 2\,{q^4} + {q^5} \right) \,{\bd 4} + \cr &&
\left( q + 2\,{q^2} + 2\,{q^3} + {q^4} \right) \,{\bd {16}} + \cr &&
\left( q + {q^2} + {q^3} \right) \,{\bd {20}} + 
\left( q + {q^2} \right) \,{\bd {40}} + {\bd {64}}, \cr
[{\bd 4}^2 \otimes {\bd 5}^2]&=&
\left( {q^3} + {q^5} \right) \,{\bd 1} + 
\left( 2\,{q^3} + {q^4} + {q^5} \right) \,{\bd 5} + \cr &&
\left( 2\,{q^2} + {q^3} + 2\,{q^4} \right) \,{\bd {10}} + 
\left( q + {q^2} + {q^3} \right) \,{\bd {14}} + \cr &&
\left( q + 2\,{q^2} + {q^3} \right) \,{\bd {35}} + q\,{\bd {30}} + q\,{\bd {35}_{a}} + {\bd {81}},
\cr
[{\bd 4} \otimes {\bd 5}^3]&=&
\left( {q^3} + {q^4} + {q^5} + {q^6} \right) \,{\bd 4} + \cr &&
\left( {q^2} + 2\,{q^3} + 2\,{q^4} + {q^5} \right) \,{\bd {16}} + \cr &&
\left( {q^2} + {q^3} + {q^4} \right) \,{\bd {20}} + 
\left( q + {q^2} + {q^3} \right) \,{\bd {40}} + \cr &&
\left( q + {q^2} \right) \,{\bd {64}} + {\bd {80}}, \cr
[{\bd 5}^4]&=&
\left( {q^4} + {q^6} + {q^8} \right) \,{\bd 1} + {q^6}\,{\bd 5} + \cr &&
\left( {q^3} + {q^4} + 2\,{q^5} + {q^6} + {q^7} \right) \,{\bd {10}} + \cr &&
\left( {q^2} + {q^3} + 2\,{q^4} + {q^5} + {q^6} \right) \,{\bd {14}} + 
\left( {q^3} + {q^4} + {q^5} \right) \,{\bd {35}} + \cr &&
\left( {q^2} + {q^4} \right) \,{\bd {35}_{a}} + 
\left( q + {q^2} + {q^3} \right) \,{\bd {81}} + {\bd {55}}.
\end{eqnarray*}

\vfill \break

\end{document}